\documentstyle[11pt]{article}
\newcommand{\E}{{\cal E}}

\newtheorem{theorem}{Theorem}[section]

\def\whitebox{{\hbox{\hskip 1pt
 \vrule height 6pt depth 1.5pt
 \lower 1.5pt\vbox to 7.5pt{\hrule width
    3.2pt\vfill\hrule width 3.2pt}%
 \vrule height 6pt depth 1.5pt
 \hskip 1pt } }}
\def\qed{\ifhmode\allowbreak\else\nobreak\fi\hfill\quad\nobreak
     \whitebox\medbreak}

\newcommand{\ignore}[1]{}

\begin {document}

\baselineskip 16pt
\title{A note on ``Extremal graphs with bounded\\ vertex bipartiteness number"}

 \author{\small   Jia\textrm{-}Bao \ Liu$^{a}$, \ \ Shaohui \ Wang$^{b,}$\thanks{Corresponding author. \newline
  \E-mail:liujiabaoad@163.com(J.-B. Liu),  shaohuiwang@yahoo.com(S. Wang).}\\
 \small  $^{a}$ School of Mathematics and Physics,
Anhui Jianzhu University, Hefei 230601, P.R. China\\
 \small  $^{b}$ Department of Mathematics and Computer Science, Adelphi University, Garden City, NY 11530, USA\\}

\date{}
\maketitle
\begin{abstract}
 This paper is devoted to present two counterexamples to the theorem
from \cite{MK} Maria R., Katherine T. M., Bernardo S. M., Extremal
graphs with bounded vertex bipartiteness number, Linear Algebra
Appl. 493 (2016) 28-36. Moreover, the corrected theorem and proof
are presented.

\end{abstract}

\section{ Introduction}
Throughout this paper we are concerned with finite undirected
connected simple graphs. Let $G$ be a graph with vertices labelled
$1,2,\dots,n$. The adjacency matrix $A_G$ of $G$ is an $n\times n$
matrix with the $(i,j)$-entry equal to 1 if vertices $i$ and $j$
are adjacent and 0 otherwise. The spectrum of a matrix $M$,
denoted by $\sigma(M)$ is the multiplicities of the eigenvalues,
which are represented in $\sigma(M)$ as powers in square brackets,
e.g., $\sigma(M)  = \{\lambda_1^{ [m_1]}, \lambda_2^{ [m_2]},
\cdots , \lambda_q^{ [m_q]} \}$ indicates that $\lambda_1$ has
multiplicity $m_1$, $\lambda_2$ has multiplicity $m_2$, and so on.
The spectral radius of $M$ is $\rho(M) = max\{|\lambda| : \lambda
\in \sigma(M)\}.$

We use $\overline{G}$ to denote the complement of $G$. The
complete graph $K_n$ is a graph on $n$ vertices such that any two
distinct vertices are connected. Let $K_{s,t}$ denote the complete
bipartite graph whose partition classes have orders $s$ and $t.$
For two vertex-disjoint graphs $G_1$ and $G_2$, the join $G_1 \vee
G_2$ is the graph such that $V(G_1 \vee G_2) = V(G_1) \cup V(G_2)$
 and
 $E(G_1 \vee G_2) = E(G_1)\cup E(G_2) \cup \{uv\mid
u\in V(G_1), v\in V(G_2)\}.$
 For the
underlying graph theoretical definitions and notations we follow
{\rm \cite{BM,MK}.

 The fewest number of vertices whose deletion
yields a bipartite graph from $G$ was defined by Fallat and
Fan to be the vertex bipartiteness of $G$ and it is
denoted by $v_b(G).$

Let $k$  be a natural number such that $3 \leq  n-k.$ The set
$\Sigma_k(n)$ is defined by

 $\{G = (V (G),E (G)) : G $ is connected,
$\mid V (G)\mid = n$ and $v_b(G)\leq k\}.$

 The authors of
\cite{MK}
 identified the graph in $\Sigma_k(n)$ with maximum spectral radius and
maximum signless Laplacian spectral radius. They obtained the
following theorem.

%
%

\begin{theorem} {\rm (Theorem 3. in
\cite{MK})} Let  $1 \leq k\leq n-3.$ Then the following hold.

(a) If $n-k$ is even, then $\rho_{\widehat{G}}\geq \rho_G$ holds
for all graphs $G\in \Sigma_k(n)$, where
 $$  \widehat{G} =K_{k}\vee\big(\overline{K_{\frac{n-k}{2}}}\vee
\overline{K_{\frac{n-k}{2}}}\big) \in \Sigma_k(n).$$ Equality
holds if and only if $ G = \widehat{G}.$ The expression for
$\rho_{\widehat{G}}\geq \rho_G$ is given by
$$\rho(\widehat{G})=\frac{1}{2}\Bigg(\frac{n+k-2}{2}+\sqrt{k^2+(3k-1)(n-k)+\Big(\frac{n-k+2}{2}\Big)^2} ~\Bigg).$$

(b) If $n-k$ is odd, then $\rho(\widehat{G})\geq \rho(G)$ holds
for all graphs $G\in \Sigma_k(n)$, where
 $$  \widehat{G} =K_{m}\vee\big(\overline{K_{\lfloor\frac{n-k}{2}\rfloor}}\vee
\overline{K_{\lfloor\frac{n-k+1}{2}\rfloor}}\big) \in
\Sigma_k(n).$$ Equality holds if and only if $ G = \widehat{G}.$
The eigenvalue $\rho(\widehat{G})$ corresponds to the maximal
eigenvalue of $\widetilde{\Upsilon}$ in
 $$ \widetilde{\Upsilon}= \left(%
\begin{array}{ccc}
  k-1 & s+1 &  s \\
 k &  0 &  s \\
    k & s+1  & 0 \\
\end{array}%
\right).$$
\end{theorem}

 However, for Theorem 1.1, we have the
following counterexamples.

\section{Some counterexamples}

\noindent In this section, we propose two examples to show the
bound of spectral radius is incorrect in Theorem 3 (a) of
\cite{MK}.

\vspace{5pt} \noindent {\bf  Example 2.1}~
 Given a graph $  \widehat{G} =K_{4}\vee\big(\overline{K_{3}}\vee
\overline{K_{3}}\big),$ by routine calculations,  we can obtain
$$\sigma(A_{K_{4}\vee\big(\overline{K_{3}}\vee
\overline{K_{3}}\big)})  = \Big\{(-3)^{ [1]}, (-1)^{ [3]},  0^{
[4]}, (3-2\sqrt{6}~)^{ [1]}, (3+2\sqrt{6}~)^{ [1]} \Big\}.$$ Then
the spectral radius $\rho(A_{K_{4}\vee\big(\overline{K_{3}}\vee
\overline{K_{3}}\big)})$ is $3+2\sqrt{6}$.

Note that $n =10, k=4$. Based on Theorem 1.1, one can get
\begin{eqnarray} \nonumber
\rho(A_{K_{4}\vee\big(\overline{K_{3}}\vee \overline{K_{3}}\big)})
&=&
\frac{1}{2}\Bigg(\frac{n+k-2}{2}+\sqrt{k^2+(3k-1)(n-k)+\Big(\frac{n-k+2}{2}\Big)^2}
~\Bigg)\\  \nonumber
 &=&
3+ \frac{7\sqrt{2}}{2}. \nonumber\end{eqnarray} \nonumber
 It is a contradiction.

\vspace{5pt} \noindent {\bf  Example 2.2}~
 Given a graph $  \widehat{G} =K_{5}\vee\big(\overline{K_{3}}\vee
\overline{K_{3}}\big).$ Similarly, we can obtain
$$\sigma(A_{K_{5}\vee\big(\overline{K_{3}}\vee
\overline{K_{3}}\big)})  = \Big\{(-3)^{ [1]}, (-2)^{ [1]}, (-1)^{
[4]}, 0^{ [4]},    9^{ [1]} \Big\}.$$
Then the spectral radius
$\rho(A_{K_{5}\vee\big(\overline{K_{3}}\vee
\overline{K_{3}}\big)})$ is
 $9$.

Note that $n =11, k=5$. According to Theorem 1.1, one can arrive
at
\begin{eqnarray} \nonumber \rho(A_{K_{5}\vee\big(\overline{K_{3}}\vee
\overline{K_{3}}\big)}) &=&
\frac{1}{2}\Bigg(\frac{n+k-2}{2}+\sqrt{k^2+(3k-1)(n-k)+\Big(\frac{n-k+2}{2}\Big)^2}
~\Bigg)\\ \nonumber &=& \frac{7}{2}+ \frac{5\sqrt{5}}{2},
\nonumber
\end{eqnarray} \nonumber
 which also arrives at a
contradiction.

We present the corresponding theorem in next section, which
corrects Theorem 1.1.

\section{Main results}
In this section, we first recall the known property of quotient
matrix that will be used later.

\begin{theorem} (See~\cite{WH}.) Suppose $\overline{M}$ is the quotient matrix of a
partitioned symmetric matrix $M$, then the eigenvalues of
$\overline{M}$ interlace the eigenvalues of $M.$ Moreover, if the
interlacing is tight, then the partition of $M$ is regular. On the
other hand, if the $M$ is regularly partitioned, then the
eigenvalues of $\overline{M}$ are eigenvalues of $M$.
\end{theorem}

In what follows, we compute the eigenvalues of $\widehat{G}$ and
obtain the following result.

\begin{theorem}  Let  $1 \leq k\leq n-3.$ Then the following hold.

(a) If $n-k$ is even, then $\rho(A_{\widehat{G}})\geq \rho(A_G)$
holds for all graphs $G\in \Sigma_k(n)$, where
 $$  \widehat{G} =K_{k}\vee\big(\overline{K_{\frac{n-k}{2}}}\vee
\overline{K_{\frac{n-k}{2}}}\big) \in \Sigma_k(n).$$ Equality
holds if and only if $ G \cong \widehat{G}.$ The expression for
$\rho(A_{\widehat{G}})$ is given by
$$\rho(A_{\widehat{G}})=\frac{n+k-2+\sqrt{-7k^2+10kn-12k+n^2+4n+4}}{4}.$$

(b) If $n-k$ is odd, then $\rho(A_{\widehat{G}})\geq \rho(A_G)$
holds for all graphs $G\in \Sigma_k(n)$, where
 $$  \widehat{G} =K_{k}\vee\big(\overline{K_{\lfloor\frac{n-k}{2}\rfloor}}\vee
\overline{K_{\lfloor\frac{n-k+1}{2}\rfloor}}\big) \in
\Sigma_k(n).$$ Equality holds if and only if $ G \cong
\widehat{G}.$ The eigenvalue $\rho(A_{\widehat{G}})$ corresponds
to the maximal eigenvalue of $\widetilde{\Upsilon}$ in
 $$ \widetilde{\Upsilon}= \left(%
\begin{array}{ccc}
  k-1 & s+1 &  s \\
 k &  0 &  s \\
    k &  s+1  & 0 \\
\end{array}%
\right).$$
\end{theorem}

\noindent{\bf Proof.} By a suitable labelling, for $ \widehat{G}
=K_{k}\vee\big(\overline{K_{\frac{n-k}{2}}}\vee
\overline{K_{\frac{n-k}{2}}}\big)$, the matrix $A_{\widehat{G}}$
 has regular partitioning into blocks.

(a) If $n-k$ is even, $r=s =\frac{n-k}{2},$ then by Theorem 3.1,
$\sigma(\widetilde{\Omega})\subseteq \sigma(A_{\widehat{G}})$,
 where
  $$ \widetilde{\Upsilon}= \left(%
\begin{array}{ccc}
  k-1 & r &  s \\
 k &  0 &  s \\
    k &  r  & 0 \\
\end{array}%
\right).$$

Note that an irreducible nonnegative matrix, the only eigenvalue
with a nonnegative eigenvector is its maximal eigenvalue
(see~\cite{HM}). If $M$ is an irreducible nonnegative matrix then
the quotient matrix $\overline{M}$ is also irreducible and
nonnegative. Consequently, the maximal eigenvalue of
$\overline{M}$ also is the maximal eigenvalue of $M$.

We now compute the eigenvalues of $\widetilde{\Upsilon}$.

\begin{eqnarray} \nonumber
 \left| \widetilde{\Upsilon}-\lambda\cdot I
\right| &=&  \left| %
\begin{array}{ccc}
  k-1-\lambda & r &  s \\
 k &  -\lambda &  s \\
    k &  r  & -\lambda \\
\end{array}%
\right|\\
 \nonumber &= &  \left| %
\begin{array}{ccc}
  k-1-\lambda & \frac{n-k}{2} &  \frac{n-k}{2} \\
 k &  -\lambda &  \frac{n-k}{2} \\
    k &  \frac{n-k}{2}  & -\lambda \\
\end{array}%
\right|\\
 \nonumber
 &= & (\lambda+\frac{n-k}{2})\Big[ \lambda^2+(1-\frac{n-k}{2}-k)\lambda-(k+1)\frac{n-k}{2}\Big] .\\
 \nonumber
\end{eqnarray}

We can obtain $$\lambda_1=\frac{k-n}{2},
\lambda_{2,3}=\frac{n+k-2\pm\sqrt{-7k^2+10kn-12k+n^2+4n+4}}{4}.$$
$$\sigma (\widetilde{\Upsilon})=\Bigg\{\frac{k-n}{2}, \frac{n+k-2\pm\sqrt{-7k^2+10kn-12k+n^2+4n+4}}{4}\Bigg\}.$$
Hence, the spectral radius $\rho(A_{\widehat{G}})$ is the maximal
eigenvalue of $\sigma (\widetilde{\Upsilon})$.
The expression for
$\rho(A_{\widehat{G}})$ is given by
$$\rho(\widehat{G})=\frac{n+k-2+\sqrt{-7k^2+10kn-12k+n^2+4n+4}}{4}.$$

 For the proof of (b), readers can
refer to the proof of Theorem 3 (b)~\cite{MK}, which is correct.
Hence, we omit it here. \qed

 \vspace{5pt} \noindent
{\bf Acknowledgments}\

The work is partly supported by the National Science Foundation of
China under Grant No. 11601006.

\end{document}